\documentclass[12pt]{article}
\usepackage{amsfonts}
\usepackage{amsthm,amsmath,amssymb,amsfonts}\newtheorem{theorem}{Theorem}[section]

\newtheorem{lemma}[theorem]{Lemma}

\newtheorem{assumption}[theorem]{Assumption}
\newtheorem{remark}[theorem]{Remark}

\begin{document}

\title{Simultaneous Identification of the Diffusion Coefficient and the 
Potential for the Schr\"{o}dinger Operator with only one Observation}

\author{Laure Cardoulis \thanks{Universit\'e de Toulouse, UT1 CEREMATH, CNRS, Institut de Math\'ematiques de Toulouse, UMR 5219, 21 Allées de Brienne 31042 Toulouse, France} and Patricia Gaitan \thanks{ Universit\'e d'Aix-Marseille, IUT Aix-en-Provence Avenue Gaston Berger 413 av Gaston Berger, 13625 Aix-en-Provence
et LATP, UMR CNRS 6632, 39, rue Joliot Curie, 13453 Marseille Cedex 13, France}}
\maketitle
\begin{abstract}
This article is devoted to prove a stability result for two
independent coefficients for a 
Schr\"{o}dinger operator in an unbounded strip. The result is obtained with only one observation on 
an unbounded subset of the boundary and the data of the solution at a fixed time on the whole domain.
\end{abstract}

\section{Introduction}

\noindent

Let $\Omega= \mathbb{R} \times (d,2d)$ be an unbounded strip of
$\mathbb{R}^2$ with a fixed width $d>0$.
Let $\nu$ be the outward unit normal to $\Omega$ on $\Gamma=\partial \Omega$. 
We denote $x=(x_1,x_2)$ and $\Gamma=\Gamma^+ \cup \Gamma^-$, where
$\Gamma^+=\{x \in \Gamma ; \; x_2=2d\}$ and 
$\Gamma^-=\{x \in \Gamma ; \; x_2=d\}$. 
We consider the following Schr\"{o}dinger equation
\begin{equation}\label{systq}
\left \{ \begin{array}{ll}
Hq:=i \partial_t q +a \Delta q+ b q=0 \;\mbox{ in }\; \Omega \times (0,T),\\
q(x,t)=F(x,t) \;\mbox{ on }\; \partial\Omega \times (0,T),\\
q(x,0)=q_0(x) \;\mbox{ in }\; \Omega,
\end{array}\right.
\end{equation}
where $a$ and $b$ are real-valued functions such that 
$a \in \mathcal{C}^3(\overline{\Omega})$, 
$b \in \mathcal{C}^2(\overline{\Omega})$ and 
$a(x) \ge a_{min}>0$. Moreover,
we assume that $a$ is bounded and $b$ and all its derivatives up to order two
are bounded.
If we assume that $q_0$ belongs to $H^4(\Omega)$ and $F \in H^2(0,T,H^{2}
(\partial \Omega)) \cap  H^1(0,T,H^{4}
(\partial \Omega)) \cap  H^3(0,T,L^{2}
(\partial \Omega))$, then (\ref{systq}) admits a solution
in $H^1(0,T, H^{2}(\Omega))\cap  H^2(0,T,L^{2}(\Omega))$.\\
Our problem can be stated as follows:\\  \noindent
Is it possible to determine the coefficients $a$ and $b$ from the measurement of
$\partial_{\nu}(\partial_t^2 q)$ on $\Gamma^+$?\\ \\ \noindent
Let $q$ (resp. $\widetilde{q}$) be a solution of (\ref{systq}) associated with
($a$, $b$, $F$, $q_0$) (resp. ($\widetilde{a}$, $\widetilde{b}$, $F$, $q_0$)).
We assume that 
$q_0$ is a real valued function.

\noindent
Our main result is 
\begin{eqnarray*}
\|a-\widetilde{a}\|^2_{L^2(\Omega)}+ \|b-\widetilde{b}\|^2_{L^2(\Omega)} & \leq & C\|\partial_{\nu}(\partial_t^2 q)-
\partial_{\nu}(\partial_t^2 \widetilde{q})\|^2_{L^2((-T,T)\times \Gamma^+)}\\
& + & C\sum_{i=0}^2\|\partial_t^i(q-\widetilde{q})(\cdot,0)\|^2
_{H^2(\Omega)},
\end{eqnarray*}
where $C$ is a positive constant which depends on $(\Omega, \Gamma, T)$
and where the above norms are weighted Sobolev norms. \\
 \noindent
This paper is an improvement of the work \cite{LMP} in the sense that we 
simultaneously determine with only one observation, two independent coefficients, 
the diffusion coefficient and the potential. We use for that two important tools: 
Carleman estimate (\ref{carl-est-2}) and Lemma \ref{lemBK}.\\
Carleman inequalities constitute a very efficient tool to derive 
observability estimates.
The method of Carleman estimates  
has been introduced in the field of inverse problems by 
Bukhgeim and Klibanov (see \cite{B}, \cite{BK}, \cite{K1}, \cite{K2}). Carleman estimates techniques are presented in \cite{KT} for standard coefficients 
inverse problems for both linear and non-linear partial differential 
equations. These methods give a local Lipschitz stability around a single 
known solution. \\
A lot of works using the same strategy concern the wave 
equation (see \cite{LTY}, \cite{Be}, \cite{BMO}) and the heat equation 
(see \cite{PY}, \cite{IY01}, \cite{BGLR}). For the determination of a 
time-independent potential in Schr\"{o}dinger evolution equation, 
we can refer to \cite{BP} for bounded domains and \cite{LMP} for unbounded
domains. We can also cite \cite{MOR} where the authors use weight functions 
satisfying a relaxed pseudo-convexity condition which allows to prove 
Carleman inequalities with less restrictive boundary observations.\\ 
Up to our knowledge, there are few results concerning the 
simultaneous identification of two coefficients with only one
observation. In \cite{LP} a stability result is 
given for the particular case where each coefficient only depends on one 
variable ($a=a(x_2)$ and $b=b(x_1)$) for the operator 
$i\partial _t q+ \nabla \cdot (a\nabla q) +bq$ in an unbounded strip of 
$\mathbb{R}^2.$ The authors give a stability result for 
the diffusion coefficient $a$ and the potential $b$ with only one observation 
in an unbounded part of the boundary.\\ 
 \noindent
A physical background could be the reconstruction of the diffusion coefficient 
and the potential in a strip in geophysics. There are also applications in quantum mechanics: 
inverse problems associated with curved quantum guides (see \cite{CDFK}, \cite{DE}, \cite{DEK}).\\
\noindent
This paper is organized as follows. Section $2$ is devoted to some usefull estimates.
We first give an adapted global Carleman estimate for the operator $H$.
We then recall the crucial Lemma given in \cite{KT}.
In Section $3$ we state and prove our main result.

\section{Some Usefull Estimates}
\subsection{Global Carleman Inequality}
%
Let $a$ be a real-valued function in $\mathcal{C}^3(\overline{\Omega})$
and  $b$ be a real-valued function in $\mathcal{C}^2(\overline{\Omega})$ 
such that
\begin{assumption} 
\label{ab}
\begin{itemize}
\item $a \ge a_{min}>0$, $a$ and all its derivatives up to order three are 
bounded,
\item $b$ and its derivatives up to order two are bounded.
\end{itemize}
\end{assumption}
Let $q(x,t)$ be a function equals to zero on $\partial \Omega \times (-T,T)$ and solution
of the Schr\"{o}dinger equation 
$$i \partial_t q +a \Delta q+ b q=f.$$
We prove here a global Carleman-type estimate for $q$
with a single observation 
acting on a part $\Gamma^+$ of the boundary $\Gamma$ 
in the right-hand side of the estimate.\\
Note that this estimate is quite similar to the one obtained in \cite{LMP}, but the 
computations are different. Indeed, the weigth function $\beta$ does not satisfy the 
same pseudo-convexity assumptions (see Assumption \ref{funct-beta}) and the decomposition
of the operator $H$ is different (see (\ref{M1})).\\
Let $\widetilde{\beta}$ be a 
$\mathcal{C}^4(\overline{\Omega})$ positive function such that there exists positive
constants $C_0,C_{pc}$ which satisfy 
\begin{assumption}
\label{funct-beta}
\begin{itemize}
\item
$|\nabla \widetilde{\beta}| \geq C_0>0 \;\;\mbox{ in } \;\;\Omega, \;\;
{\partial}_{\nu} {\widetilde{\beta}}\leq 0\;\;\mbox{on}\;\;\Gamma^-$,
\item
$ \widetilde{\beta}$ and all its derivatives up to order four are bounded in
$\overline{\Omega}$,
\item
$2 \Re (D^2 \widetilde{\beta}(\zeta, \bar{\zeta})) 
-\nabla a \cdot \nabla \widetilde{\beta} |\zeta|^2+2a^2|\nabla \widetilde{\beta} \cdot \zeta|^2
 \geq C_{pc} |\zeta|^2$, for all $\zeta \in \mathbb{C}$
\end{itemize}
\end{assumption}
where 
$$D^2 \widetilde{\beta}=\left( \begin{array}{cc}
\partial_{x_1}(a^2\partial_{x_1} \widetilde{\beta}) & \partial_{x_1}(a^2\partial_{x_2} \widetilde{\beta})\\
\partial_{x_2}(a^2\partial_{x_1} \widetilde{\beta}) & \partial_{x_2}(a^2\partial_{x_2} \widetilde{\beta})
\end{array} \right).$$
Note that the last assertion of Assumption \ref{funct-beta} expresses
the pseudo-convexity condition for the function $\beta$.
This Assumption imposes restrictive conditions for the choice of the diffusion coefficient
$a$ in connection with the function $\widetilde{\beta}$ as in \cite{LMP}. \\
Note that there exist functions satisfying such assumptions. Indeed  
if we assume that $\tilde{\beta}(x):=\tilde{\beta}(x_2)$, these conditions 
can be written in the following form:
$$A=2\partial_{x_2} (a^2 \partial_{x_2} \tilde{\beta}) -\partial_{x_2}a \ 
\partial_{x_2} \tilde{\beta} +2a^2 (\partial_{x_2} \tilde{\beta})^2 \geq cst 
>0$$ 
and 
$$-\frac{(\partial_{x_1}(a^2 \partial_{x_2} \tilde{\beta}))^2}{A} -
\partial_{x_2}a \ \partial_{x_2} \tilde{\beta} \geq cst>0.$$
 For example $\tilde{\beta}(x)=
e^{-x_2}$ with $a(x)=\frac{1}{2} (x_2 ^2 +5)$ satisfy the previous conditions 
(with $x_2 \in (d,2d)$).

\noindent
Then, we define $\beta= \widetilde{\beta}+K$ with
$K= m \|\widetilde{\beta}\|_{\infty}$ and $m>1$. For $\lambda> 0$ and $t \in
(-T,T)$, we define the following weight functions
$$
  \varphi(x,t)=\frac{e^{\lambda \beta(x)}}{(T+t)(T-t)},
  \quad \quad \eta(x,t)=\frac{e^{2\lambda K} -e^{\lambda
  \beta(x)}}{(T+t)(T-t)}.
$$
We set $\psi=e^{-s \eta}q$, $M \psi = e^{-s \eta} H(e^{s \eta} \psi)$ for $s>0$. Let $H$ be the operator defined by
\begin{equation} \label{H}
Hq:=i\partial_t q + a \Delta q + b q  \;\mbox{ in }\; \Omega \times (-T,T).
\end{equation}
Following \cite{BP}, we introduce the operators  :
\begin{eqnarray} \label{M1}
  M_1\psi  & : = & i\partial_t \psi+ a \Delta \psi+s^2 a |\nabla \eta |^2 \psi + (b-s\nabla \eta \cdot \nabla a) \psi,\\ \nonumber
  M_2\psi  & : = & is \partial_t\eta \psi+2 a s \nabla \eta  \cdot \nabla \psi +s \nabla \cdot (a \nabla \eta) \psi.
\end{eqnarray}
Then
\begin{eqnarray*}
\int_{-T}^T  \int_{\Omega}|M\psi|^2dx  \ dt & = & 
\int_{-T}^T   \int_{\Omega}|M_1\psi|^2dx \ dt
+\int_{-T}^T   \int_{\Omega}|M_2\psi|^2dx  \ dt\\
 & + &  2 \Re (\int_{-T}^T   \int_{\Omega}M_1\psi\  \overline{M_2\psi}\ dx  \ dt),
\end{eqnarray*}
where $\overline{z}$ is the conjugate of $z$, $\Re \ (z)$ its real part
and $\Im \ (z)$ its imaginary part.
Then the following result holds.\\
\begin{theorem}
\label{th-Carl} 
Let $H$, $M_1$, $M_2$ be the operators defined respectively by
(\ref{H}), (\ref{M1}). We assume that Assumptions \ref{ab} and \ref{funct-beta}
are satisfied. 
Then there exist $\lambda_0> 0$, $s_0>0$ and a positive
constant $C=C(\Omega, \Gamma,T)$ such that, for any $\lambda \ge
\lambda_0$ and any $s \ge s_0 $, the next inequality holds:
$$
s^{3} \lambda^{4}\int_{-T}^T   \int_{\Omega} e^{-2s \eta}  |q|^2\ d x  \ d t  
+s \lambda \int_{-T}^T   \int_{\Omega} e^{-2s \eta}  |\nabla q|^2\ dx \ d t
  +  \|M_1(e^{-s \eta}q)\|^2_{L^2(\Omega \times (-T,T))}$$
\begin{eqnarray}\label{carl-est-1}
+ \|M_2(e^{-s \eta}q)\|^2_{L^2(\Omega \times (-T,T))}
  \leq C 
   s \lambda \int_{-T}^T \int_{\Gamma^+} e^{-2s \eta} |\partial_{\nu} q|^2\ 
\partial_{\nu} \beta \ d \sigma \ d t \\ \nonumber
  +\int_{-T}^T   \int_{\Omega} e^{-2s \eta}\ | H q |^2\ d x \ d t,
  \end{eqnarray}
for all $q$ satisfying $q\in
L^2(-T,T;H^1_0(\Omega)\cap H^2(\Omega))\cap H^1(-T,T;L^2(\Omega)),$ 
$\partial_{\nu} q\in L^2(-T,T;L^2(\Gamma)).$
Moreover we have

$$s^{3} \lambda^{4}\int_{-T}^T   \int_{\Omega} e^{-2s \eta}  |q|^2\ d x  \ d t  
+s \lambda \int_{-T}^T   \int_{\Omega} e^{-2s \eta}  |\nabla q|^2\ dx \ d t 
+  \|M_1(e^{-s \eta}q)\|^2_{L^2(\Omega \times (-T,T))}$$
\begin{eqnarray}\label{carl-est-2} 
+ 
\|M_2(e^{-s \eta}q)\|^2_{L^2(\Omega \times (-T,T))} +s^{-1} \lambda^{-1} 
\int_{-T}^T   \int_{\Omega} e^{-2s \eta}  |i\partial_t q+a\Delta q|^2\ d x  \ d t\\ \nonumber
  \leq C \left[
   s \lambda \int_{-T}^T \int_{\Gamma^+} e^{-2s \eta} |\partial_{\nu} q|^2\ \partial_{\nu} \beta \ d \sigma \ d t
  +\int_{-T}^T   \int_{\Omega} e^{-2s \eta}\ | H q |^2\ d x \ d t\right].
  \end{eqnarray}
\end{theorem}
{\bf Proof:}\\
We have to estimate the scalar product
$$\Re \left(\int_{-T}^T   \int_{\Omega}M_1\psi\  \overline{M_2\psi}\ dx  \ dt\right)=
\sum_{i=1}^4 \sum_{j=1}^3 I_{ij}$$
with
{\scriptsize{$$I_{11}=\Re \left( \int_{-T}^T  \int_{\Omega} (i\partial_t \psi)
(-is\partial_t\eta \ \overline{\psi})\  dx  \ dt \right),\ \ \ 
I_{12}=\Re \left( \int_{-T}^T  \int_{\Omega} (i\partial_t \psi)
(2as\nabla\eta \cdot\nabla \overline{\psi})\  dx  \ dt \right),$$
$$I_{13}=\Re \left( \int_{-T}^T  \int_{\Omega} (i\partial_t \psi)
(s \nabla\cdot(a \nabla \eta) \overline{\psi})\  dx  \ dt \right),\ \ \ 
I_{21}=\Re \left( \int_{-T}^T  \int_{\Omega} (a \Delta \psi)
(-is\partial_t\eta \ \overline{\psi})\  dx  \ dt \right),$$
$$I_{22}=\Re \left( \int_{-T}^T  \int_{\Omega} (a \Delta \psi)
(2 a s \nabla\eta \cdot \nabla \overline{\psi})\  dx  \ dt \right),\ \ \ 
I_{23}=\Re \left( \int_{-T}^T  \int_{\Omega} (a \Delta\psi)
(s\nabla \cdot(a \nabla\eta) \ \overline{\psi})\  dx  \ dt \right),$$
$$I_{31}=\Re \left( \int_{-T}^T  \int_{\Omega} (s^2 a |\nabla \eta|^2 \psi)
(-is\partial_t\eta \ \overline{\psi})\  dx  \ dt \right),\ \ \    
I_{32}=\Re \left( \int_{-T}^T  \int_{\Omega} (s^2 a |\nabla \eta|^2 \psi)
(2as\nabla\eta \cdot \nabla \overline{\psi})\  dx  \ dt \right),$$
$$I_{33}=\Re \left( \int_{-T}^T  \int_{\Omega} (s^2 a |\nabla \eta|^2 \psi)
(s\nabla \cdot(a \nabla\eta) \ \overline{\psi})\  dx  \ dt \right),\ \ \ 
I_{41}=\Re \left( \int_{-T}^T  \int_{\Omega} ((b-s\nabla \eta \cdot \nabla a)\psi)
(-is\partial_t\eta \ \overline{\psi})\  dx  \ dt \right),$$
$$ I_{42}=\Re \left( \int_{-T}^T  \int_{\Omega} ((b-s\nabla \eta \cdot \nabla a)\psi)
(2 a s \nabla\eta \cdot \nabla \overline{\psi})\  dx  \ dt \right),\ \ \ 
I_{43}=\Re \left( \int_{-T}^T  \int_{\Omega} ((b-s\nabla \eta \cdot \nabla a)\psi)
(s\nabla \cdot(a \nabla\eta) \ \overline{\psi})\  dx  \ dt \right).$$}}

\noindent
Following \cite{BP}, using integrations by part and Young estimates, we get (\ref{carl-est-1}). Moreover from (\ref{M1}) we 
have:
$$i\partial_t q+a \Delta q=M_1 q-s^2 a|\nabla  \eta|^2 q+(b-s \nabla \eta 
\cdot \nabla a)q.$$
So 
\begin{eqnarray*}
i\partial_t q+a \Delta q & = & e^{s\eta} M_1(e^{-s\eta} q)+is \partial_t \eta q 
-ae^{s\eta} \Delta(e^{-s\eta})q -2a e^{s\eta} \nabla(e^{-s\eta}) 
\cdot \nabla q \\
&-& s^2 a|\nabla  \eta|^2 q+(b-s \nabla \eta 
\cdot \nabla a)q.
\end{eqnarray*}
And we deduce (\ref{carl-est-2}) from (\ref{carl-est-1}).
\begin{flushright}
\rule{.05in}{.05in}
\end{flushright}
%
\subsection{The Crucial Lemma}
%
We recall in this section the proof of a very important lemma proved by Klibanov and Timonov (see for example \cite{K2}, \cite{KT}).
\begin{lemma}
\label{lemBK} There exists a positive constant $\kappa$ such that 
\[
\int_{-T}^{T} \int_{\Omega}\left\vert \int^t_{0} q(x,\xi) d\xi\right\vert^2 e^{-2s\eta} dxdt
\le \frac{\kappa}{s}\int_{-T}^{T} \int_{\Omega } \vert q(x,t)\vert^2 e^{-2s\eta} dxdt, 
\]
for all $s > 0$.
\end{lemma}

{\bf Proof :}\\
By the Cauchy-Schwartz inequality, we have
\begin{equation} \label{estBK}
\int_{-T}^{T} \int_{\Omega} \left\vert \int^t_{0} q(x,\xi) d\xi\right\vert^2 e^{-2s\eta} dxdt
\le \int_{\Omega} \int^T_{-T}  \vert t\vert
\left\vert \int^t_{0} \vert q(x,\xi) \vert^2 d\xi\right\vert e^{-2s\eta}
dxdt 
\end{equation}
$$\leq \int_{\Omega} \int^T_{0} t
\left( \int^t_{0} \vert q(x,\xi) \vert^2 d\xi\right)e^{-2s\eta} dxdt 
+ \int_{\Omega} \int^{0}_{-T} (-t)
\left( \int_t^{0} \vert q(x,\xi) \vert^2 d\xi\right) e^{-2s\eta}
dxdt.
$$
Note that 
$$\partial_t(e^{-2s\eta(x,t)})=-2s(e^{2\lambda K}-e^{\lambda \beta(x)})\frac{2t}{(T^2-t^2)^2}e^{-2s\eta(x,t)}.$$
So, if we denote by $\alpha(x)=e^{2\lambda K}-e^{\lambda \beta(x)}$, we have
$$
t e^{-2s\eta(x,t)} = -\frac{(T^2-t^2)^2}{4 s \alpha(x)}
\partial_t(e^{-2s\eta(x,t)}).
$$
For the first integral of the right hand side of (\ref{estBK}), by integration by parts we have
$$
\int_{\Omega} \int^T_{0} t
\left( \int^t_{0} \vert q(x,\xi) \vert^2 d\xi\right)e^{-2s\eta} dxdt
= \int_{\Omega}  \int^T_{0} 
\left( \int^t_{0} \vert q(x,\xi) \vert^2 d\xi\right)
\frac{(T^2-t^2)^2 }{-4 s \alpha(x)}\partial_t(e^{-2s\eta}) dt \ dx $$
$$ = \int_{\Omega} \left[
\left( \int^t_{0} \vert q(x,\xi) \vert^2 d\xi\right)
\frac{(T^2-t^2)^2}{-4 s \alpha(x)}e^{-2s\eta} \right]^{t=T}_{t=0} dx
+ \int_{\Omega} \int^T_{0} 
\vert q(x,t) \vert^2 \frac{(T^2-t^2)^2}{4 s \alpha(x)}e^{-2s\eta} dt \ dx$$
$$+ \int_{\Omega} \int^T_{0} 
\left( \int^t_{0} \vert q(x,\xi) \vert^2 d\xi\right)
\frac{t(t^2-T^2)}{s \alpha(x)}e^{-2s\eta} dt \ dx.
$$
Here we used $\alpha(x) > 0$ for all $x \in \overline{\Omega}$ and we obtain
$$ \int_{\Omega} \int^{T}_0 t
\left( \int^{0}_t \vert q(x,\xi) \vert^2 d\xi\right)e^{-2s\eta} dxdt
\le \frac{1}{4s}\sup_{x\in \overline{\Omega}} \left(
\frac{1}{\alpha(x)}\right) \int_{\Omega} \int^{T}_0 
\vert q(x,t) \vert^2 e^{-2s\eta} (T^2-t^2)^2 dxdt. $$
Similarly for the second integral of the right hand side of (\ref{estBK}) 
$$
\int_{\Omega} \int_{-T}^0 (-t)
\left( \int^{0}_t \vert q(x,\xi) \vert^2 d\xi\right)e^{-2s\eta} dxdt
\le \frac{1}{4s}\sup_{x\in \overline{\Omega}} \left(
\frac{1}{\alpha(x)}\right) \int_{\Omega} \int_{-T}^0 
\vert q(x,t) \vert^2 e^{-2s\eta}(T^2-t^2)^2 dxdt.
$$
Thus the proof of Lemma \ref{lemBK} is completed.
\begin{flushright}
\rule{.05in}{.05in}
\end{flushright}
%
\section{Stability result}
%
In this section, we establish a stability inequality for the diffusion coefficient $a$ and the potential $b$. \\
Let $q \in \mathcal{C}^2(\Omega \times (0,T))$ be solution of 
$$
\left \{ \begin{array}{ll}
i \partial_t q +a \Delta q+ b q=0 \;\mbox{ in }\; \Omega \times (0,T),\\
q(x,t)=F(x,t) \;\mbox{ on }\; \partial\Omega \times (0,T),\\
q(x,0)=q_0(x) \;\mbox{ in }\; \Omega,
\end{array}\right.
$$
and $\widetilde{q} \in \mathcal{C}^2(\Omega \times (0,T))$ be solution of 
$$
\left \{ \begin{array}{ll}
i \partial_t \widetilde{q} +\widetilde{a} \Delta \widetilde{q}+ 
\widetilde{b} \widetilde{q}=0 \;\mbox{ in }\; \Omega \times (0,T),\\
\widetilde{q}(x,t)=F(x,t) \;\mbox{ on }\; \partial\Omega \times (0,T),\\
\widetilde{q}(x,0)=q_0(x) \;\mbox{ in }\; \Omega,
\end{array}\right.
$$
where $(a,b)$ and $(\widetilde{a},\widetilde{b})$ both satisfy Assumption 
\ref{ab}. 
\begin{assumption}\label{tildeq}
\begin{itemize}
\item All the time-derivatives up to order three and the space-derivatives up 
to order four for $\tilde{q}$ exist and are bounded.
\item There exists a positive constant $C>0$ such that $|\tilde{q}| \geq C$, 
$\displaystyle{|\partial_{t} (\frac{\Delta \tilde{q}}{\tilde{q}})|} \geq C,$ 
$|\Delta \tilde{q} | \geq C,$ 
$\displaystyle{|\partial_{t} (\frac{\tilde{q}}{\Delta \tilde{q}})|} \geq C.$
\item $q_0$ is a real-valued function. 
\end{itemize}
\end{assumption}
Since $q_0$ is a real-valued function, we can extend the function $q$ (resp. 
$\widetilde{q}$) on $\Omega \times (-T,T)$ by the formula $q(x,t)=\overline{q}
(x,-t)$ for every $(x,t) \in \Omega \times (-T,0).$ Note that this extension 
satisfies the previous Carleman estimate. 
Our main stability result is 

\begin{theorem} \label{th-stab}
 Let $q$ and $\widetilde{q}$ be solutions of (\ref{systq}) in $ \mathcal{C}^2(\Omega \times (0,T))$
such that $q-\widetilde{q} \in H^2((-T,T);H^2(\Omega))$.
We assume that Assumptions \ref{ab}, 
\ref{funct-beta}, \ref{tildeq} are satisfied.
  Then there exists a positive constant
  $C=C(\Omega, \Gamma,T)$ 
such that for $s$ and $\lambda$ large enough,
$$\int_{-T}^{T} \int_{\Omega} e^{-2s\eta} (|\widetilde{a}-a|^2  
+|\widetilde{b}-b |^2) \ dx \ dt  
\leq C  s \lambda ^2 \int_{-T}^T \int_{\Gamma^+} \varphi \ e^{-2s\eta} 
\partial_{\nu} \beta \ |\partial_{\nu} (\partial_t^2 q -
\partial_t^2 \widetilde{q})|^2 \ d\sigma \ dt 
$$
$$+ C\lambda 
\int_{-T}^{T} \int_{\Omega} e^{-2s\eta} \bigg(
\sum_{i=0}^2 |\partial_t ^i (q-\tilde{q})(.,0)|^2 + 
\nabla (q-\tilde{q})(.,0)|^2 $$
$$ + |\partial_t \nabla(q-\tilde{q})(.,0)|^2 +|\partial_t 
\Delta (q-\tilde{q})(.,0)|^2 \bigg) \ dx \ dt.$$
Therefore 
\begin{eqnarray*}
\|a-\widetilde{a}\|^2_{L^2(\Omega)}+ \|b-\widetilde{b}\|^2_{L^2(\Omega)} & \leq & C\|\partial_{\nu}(\partial_t^2 q)-
\partial_{\nu}(\partial_t^2 \widetilde{q})\|^2_{L^2((-T,T)\times \Gamma^+)}\\
& + & C\sum_{i=0}^2\|\partial_t^i(q-\widetilde{q})(\cdot,0)\|^2
_{H^2(\Omega)},
\end{eqnarray*}
where the previous norms are weighted Sobolev norms.
\end{theorem}
{\bf Proof:}\\
We denote by $u=q-\tilde{q},$ $\alpha=\tilde{a}-a$ and $\gamma=\tilde{b}-b$, so 
we get:
\begin{equation}\label{systu}
\left \{ \begin{array}{ll}
i \partial_t u +a \Delta u+ b u=\alpha \Delta \tilde{q} +
\gamma \tilde{q} \;\mbox{ in }\; \Omega \times (-T,T),\\
u(x,t)=0 \;\mbox{ on }\; \partial\Omega \times (-T,T),\\
u(x,0)=0 \;\mbox{ in }\; \Omega.
\end{array}\right.
\end{equation}
The proof will be done in two steps:
in a first step we prove an estimation for $\alpha$ and in a second step for 
$\gamma.$\\
\noindent
{\bf First step:} We set $\displaystyle{u_1=\frac{u}{\tilde{q}}}.$ Then from (\ref{systu}) $u_1$ is solution of 
$$\left \{ \begin{array}{ll}
i \partial_t u_1 +a \Delta u_1+ b u_1 +A_{11}u_1 +B_{11} \cdot \nabla u_1=
\alpha \displaystyle{\frac{\Delta \tilde{q}}{\tilde{q}}} +
\gamma  \;\mbox{ in }\; \Omega \times (-T,T),\\
u_1(x,t)=0 \;\mbox{ on }\; \partial \Omega \times (-T,T)
\end{array}\right.$$
where $\displaystyle{A_{11}=i\frac{\partial_t \tilde{q}}{\tilde{q}} +a\frac{\Delta \tilde{q}}{\tilde{q}}}$ and $\displaystyle{B_{11}=\frac{2a}{\tilde{q}} \nabla \tilde{q}}.$\\
\noindent
Then defining $u_2=\partial_t u_1$ we get that $u_2$ satisfies
$$\left \{ \begin{array}{ll}
i \partial_t u_2 +a \Delta u_2+ b u_2 +\sum_{i=1}^2 A_{i2}u_i 
+\sum_{i=1}^2 B_{i2} \cdot \nabla u_i=
\alpha \displaystyle{\partial_t (\frac{\Delta \tilde{q}}{\tilde{q}})}  \;\mbox{ in }\; 
\Omega \times (-T,T),\\
u_2(x,t)=0 \;\mbox{ on }\; \partial \Omega \times (-T,T)
\end{array}\right.$$
where $A_{12}=\partial_t A_{11}, \ A_{22}=A_{11}, \ B_{12}=\partial_t B_{11}, 
\ B_{12}=B_{11}.$\\
\noindent
Now let $u_3=\displaystyle{\frac{u_2} 
{\partial_t (\frac{\Delta \tilde{q}}{\tilde{q}})}},$ then $u_3$ is solution of
\begin{equation}\label{u3}
\left \{ \begin{array}{ll}
i \partial_t u_3 +a \Delta u_3+ b u_3 +\sum_{i=1}^3 A_{i3}u_i 
+\sum_{i=1}^3 B_{i3} \cdot \nabla u_i=
\alpha   \;\mbox{ in }\; \Omega \times (-T,T),\\
u_3(x,t)=0 \;\mbox{ on }\; \partial \Omega \times (-T,T)
\end{array}\right.
\end{equation}
where $A_{i3}$ and $B_{i3}$ are bounded functions.\\
If we denote by 
$\displaystyle{g=\partial_t (\frac{\Delta \tilde{q}}{\tilde{q}})}$, then 
$$\displaystyle{A_{13}=\frac{1}{g} A_{12}, \ A_{23}=\frac{1}{g}A_{22},\ A_{33}=\frac{1}{g} 
(i\partial_t g +\Delta g), \ B_{13}=\frac{1}{g}B_{12},\ B_{23}=\frac{1}{g}
B_{22},\ B_{33}=\frac{2a}{g}\nabla g.}$$
\noindent
At last we define $u_4=\partial_t u_3$ and $u_4$ satisfies
$$\left \{ \begin{array}{ll}
i \partial_t u_4 +a \Delta u_4+ b u_4 +\sum_{i=1}^4 A_{i4}u_i 
+\sum_{i=1}^4 B_{i4} \cdot \nabla u_i=0   \;\mbox{ in }\; 
\Omega \times (-T,T),\\
u_4(x,t)=0 \;\mbox{ on }\; \partial \Omega \times (-T,T)
\end{array}\right.$$
where $A_{i4}$ and $B_{i4}$ are still bounded functions. Note that $A_{14}=
\partial_t A_{13},\ A_{24}=\partial_t A_{23}+A_{13},\ A_{34}=\partial_t 
A_{33} +A_{23}\partial_t g +B_{23}\cdot \nabla(\partial_t g), \
A_{44}=A_{23} g+A_{33}+B_{23}\cdot \nabla g,\ B_{14}=\partial_t B_{13}, 
\ B_{24}=\partial_t B_{23} +B_{13},\ B_{34}=\partial_t B_{33}+\partial_t g 
B_{23},\ B_{44}=B_{33}+g B_{23}.$\\
\noindent
Applying the Carleman inequality (\ref{carl-est-2}) for $u_4$ we obtain 
(for $s$ and $\lambda$ sufficiently large):
\begin{equation}
\label{u_4}
s^{3} \lambda^{4}\int_{-T}^T   \int_{\Omega} e^{-2s \eta}  |u_4|^2\ d x  \ d t  
+s \lambda \int_{-T}^T   \int_{\Omega} e^{-2s \eta}  |\nabla u_4|^2\ dx \ d t 
\end{equation}
$$ + s^{-1} \lambda^{-1}
\int_{-T}^T   \int_{\Omega} e^{-2s \eta}  |i\partial_t u_4+a\Delta u_4|^2\ d x  \ d t$$
$$  \leq C \left[
   s \lambda \int_{-T}^T \int_{\Gamma^+} e^{-2s \eta} |\partial_{\nu} u_4|^2\ \partial_{\nu} \beta \ d \sigma \ d t
  +\sum_{i=1}^3 \int_{-T}^T   \int_{\Omega} e^{-2s \eta}\ (|u_i|^2+|\nabla 
u_i|^2) \ d x \ d t\right].$$
Note that $\int_{-T}^T   \int_{\Omega} e^{-2s \eta}  |u_1|^2 
\ dx \ dt =\int_{-T}^T   \int_{\Omega} e^{-2s \eta}  |
\int_0 ^t \partial_t u_1|^2 \ dx \ dt$, so
from Lemma \ref{lemBK} we get 
$$
\int_{-T}^T   \int_{\Omega} e^{-2s \eta}  |u_1|^2 \ dx \ dt \leq 
\frac{C}{s} \int_{-T}^T   \int_{\Omega} e^{-2s \eta} |u_3|^2 \ dx \ dt $$
$$\leq  \frac{C}{s^2} \int_{-T}^T   \int_{\Omega} e^{-2s \eta} |u_4|^2 \ dx 
\ dt+ \frac{C}{s} \int_{-T}^T   \int_{\Omega} e^{-2s \eta} |u_3(.,0)|^2 \ dx 
\ dt.
$$
By the same way, we have 
$$
\int_{-T}^T   \int_{\Omega} e^{-2s \eta}  |\nabla u_1|^2 \ dx \ dt \leq 
 \frac{C}{s^2} \int_{-T}^T   \int_{\Omega} e^{-2s \eta} |\nabla u_4|^2 
\ dx \ dt$$
$$+ \frac{C}{s} \int_{-T}^T   \int_{\Omega} e^{-2s \eta} 
|\nabla u_3 (.,0)|^2 \ dx \ dt 
+
C \int_{-T}^T   \int_{\Omega} e^{-2s \eta} |\nabla u_1 (.,0)|^2 \ dx \ dt.
$$
So (\ref{u_4}) becomes
\begin{equation}
\label{u4}
s^{3} \lambda^{4}\int_{-T}^T   \int_{\Omega} e^{-2s \eta}  |u_4|^2\ d x  \ d t  
+s \lambda \int_{-T}^T   \int_{\Omega} e^{-2s \eta}  |\nabla u_4|^2\ dx \ d t 
\end{equation}
$$ + s^{-1} \lambda^{-1}
\int_{-T}^T   \int_{\Omega} e^{-2s \eta}  |i\partial_t u_4+a\Delta u_4|^2\ d x  \ d t
 \leq C 
   s \lambda \int_{-T}^T \int_{\Gamma^+} e^{-2s \eta} |\partial_{\nu} u_4|^2\ \partial_{\nu} \beta \ d \sigma \ d t$$
 $$ +C \int_{-T}^T   \int_{\Omega} e^{-2s \eta} (|u_3(.,0)|^2+|\nabla u_3 (.,0)|^2 +\nabla u_1 (.,0)|^2) \ dx \ dt.$$
\noindent
Furthermore from (\ref{u3}) we have (with $C$ a positive constant)
$$
|\alpha |^2 \leq C\left(|i\partial_t u_3 +a \Delta u_3|^2 +\sum_{i=1}^3 (|u_i|^2 +|\nabla u_i|^2) \right).
$$
Therefore for $s$ sufficiently large, from Lemma \ref{lemBK}
$$
\int_{-T}^T   \int_{\Omega} e^{-2s \eta}  |\alpha|^2  \ dx \ dt \leq  
\frac{C}{s} \int_{-T}^T   \int_{\Omega} e^{-2s \eta} 
\left( |i\partial_t u_4+a\Delta u_4|^2+|u_4|^2+|\nabla u_4|^2 \right)  \ dx \ dt $$
$$+  C  \int_{-T}^T   \int_{\Omega} e^{-2s \eta} |(i\partial_t u_3 +a\Delta u_3)(0)|^2  \ dx \ dt 
 + 
C \int_{-T}^T   \int_{\Omega} e^{-2s \eta} |\nabla u_1 (.,0)|^2  \ dx \ dt $$
 $$+  C \int_{-T}^T   \int_{\Omega} e^{-2s \eta} (|u_3(.,0)|^2+|\nabla u_3 (.,0)|^2) 
 \ dx \ dt.
$$
Using (\ref{u4}) we get
\begin{eqnarray*}
\frac{1}{\lambda} \int_{-T}^T   \int_{\Omega} e^{-2s \eta}  |\alpha|^2  \ dx \ dt
& \leq & 
C s \lambda \int_{-T}^T \int_{\Gamma^+} e^{-2s \eta} |\partial_{\nu} u_4|^2\ \partial_{\nu} \beta \ d \sigma \ d t \\
& + &
\frac{C}{\lambda} \int_{-T}^T \int_{\Omega} e^{-2s \eta} |(i\partial_t u_3 
+a\Delta u_3)(.,0)|^2  \ dx \ dt \\
& + & 
C\int_{-T}^T   \int_{\Omega} e^{-2s \eta} |\nabla u_1 (.,0)|^2  \ dx \ dt \\
& + &C\int_{-T}^T   \int_{\Omega} e^{-2s \eta} (|u_3(.,0)|^2
+|\nabla u_3 (.,0)|^2)  \ dx \ dt
\end{eqnarray*}
and then 
\begin{equation}\label{alpha}
\frac{1}{\lambda} \int_{-T}^T   \int_{\Omega} e^{-2s \eta}  |\alpha|^2 \ dx \ dt\leq 
C s \lambda \int_{-T}^T \int_{\Gamma^+} e^{-2s \eta} |\partial_{\nu} u_4|^2\ \partial_{\nu} \beta \ d \sigma \ d t
\end{equation}
$$ + C \int_{-T}^T  \int_{\Omega} e^{-2s \eta} \left( \sum_{i=0}^2 |\partial_t^i  u(.,0)|^2 
+ |\nabla u (.,0)|^2+ |\partial_t \nabla u(.,0)|^2+|\partial_t \Delta u (.,0)|^2 \right)\ dx \ dt. $$
\noindent
{\bf Second step:} By the same way we obtain an estimation of $\gamma$. We set
 $$v_1=\frac{u}{\Delta \tilde{q}},\ \ \  v_2=\partial_t v_1,\ \ \ v_3=\frac{v_2}{\partial_t (\frac{\tilde{q}}{\Delta \tilde{q}})}.$$ 
Following the same methodology as in the first step, we obtain: 
\begin{equation}\label{gamma}
\frac{1}{\lambda} \int_{-T}^T   \int_{\Omega} e^{-2s \eta}  |\gamma|^2 
\ dx \ dt
\leq 
C s \lambda \int_{-T}^T \int_{\Gamma^+} e^{-2s \eta} |\partial_{\nu} u_4|^2\ \partial_{\nu} \beta \ d \sigma \ d t\end{equation}
$$ + C \int_{-T}^T  \int_{\Omega} e^{-2s \eta} \left( \sum_{i=0}^2 |\partial_t^i  u(.,0)|^2 
+ |\nabla u (.,0)|^2+ |\partial_t \nabla u(.,0)|^2+|\partial_t \Delta u (.,0)|^2 \right) 
\ dx \ dt. $$
From (\ref{alpha}) and (\ref{gamma}) we can conclude.
\begin{flushright}
\rule{.05in}{.05in}
\end{flushright}

\begin{remark}
\begin{enumerate}
\item Note that the following function $\tilde{q} (x,t)=e^{-it}+x_2^2+5$ with 
$\displaystyle{\tilde{a}(x)=
\frac{x_2^2 +5}{2}},$ $\tilde{b}(x)=-1$ satisfies Assumption 
\ref{tildeq}. 
\item This method works for the Schr\"{o}dinger operator in the divergential 
form: $$i\partial_t q+\nabla \cdot (a\nabla q)+bq.$$ We still obtain a similar 
stability result but with more restrictive hypotheses on the regularity of 
the function $\widetilde{q}$.  
\end{enumerate}
\end{remark}

\noindent 
{\bf Acknowledgment}: We dedicate this paper to the memory of our friend and colleague Pierre Duclos, 
Professor at the University of Toulon in France.

\section*{References}



\begin{thebibliography}{99}
%
 \bibitem{BP}{\sc L. Baudouin and J.P. Puel}, 2002, { Uniqueness and stability in an inverse
problem for the Schr\"{o}dinger equation}, {\em Inverse Problems}, 18,
 1537--1554.
%
\bibitem{BMO}{\sc L. Baudouin, A. Mercado and A. Osses}, 2007, {A global Carleman estimate in a transmission wave equation and application to a 
one-measurement inverse problem}, {\em Inverse Problems}, 23,
 257--278.
%
 \bibitem{Be}{\sc M. Bellassoued}, 2004, { Uniqueness and stability in determining the speed of propagation of second order hyperbolic equation with variable 
coefficient}, {\em Appl. An.}, 83,
 983--1014.
%
 \bibitem{BGLR}{\sc A. Benabdallah, P. Gaitan and J. Le Rousseau}, 2007,  
{ Stability of discontinuous diffusion coefficients and initial conditions 
in an inverse problem for the heat equation}, {\em SIAM J. Control Optim.}, 46,
 1849--1881.
%
\bibitem{B}{\sc A.L. Bukhgeim}, 1999, Volterra Equations and Inverse Problems, {\em Inverse and Ill-Posed Problems
Series, VSP, Utrecht}.
%
\bibitem{BK}{\sc A.L. Bukhgeim and M.V. Klibanov}, 1981, Uniqueness in the large of a class of multidimensional
inverse problems, {\em Soviet Math. Dokl.}, 17, 244-247.
%
\bibitem{CDFK}{\sc B. Chenaud, P. Duclos, P. Freitas and D. Krejcirik}, {\em Geometrically induced discrete spectrum in curved tubes}, 
Diff. Geom. Appl. 23 (2005), no2, 95-105.
%
\bibitem{DE}{\sc P. Duclos and P. Exner}, {\em Curvature-Induced Bound 
States in Quantum Waveguides in Two and Three Dimensions}, Rev. Math. Phys. 
7, (1995), 73-102.
%
\bibitem{DEK}{\sc P. Duclos, P. Exner and D. Krejcirik}, {\em Bound 
States in Curved Quantum Layers}, Comm. Math. Phys. 223 (2001), 13-28.
%
\bibitem{LMP}{\sc L. Cardoulis, M. Cristofol and P. Gaitan},  2008, 
{ Inverse problem for the Schr\"{o}dinger operator in an unbounded strip}, 
{\em J. Inverse and  Ill-posed Problems}, 16, no 2, 127--146.
%
\bibitem{LP}{\sc L. Cardoulis and P. Gaitan}, 2009,  
{  Identification of two independent coefficients with one observation  
for the Schr\"{o}dinger operator in an unbounded strip}, to appear 
in {\em Comptes Rendus Acad\'emie des Sciences}.
%
\bibitem{IY01}{\sc O. Yu. Immanuvilov and M. Yamamoto}, 2001,  
{Global uniqueness and stability in determining coefficients of wave 
equations}, 
{\em Comm. Partial Diff. Equat.}, 26, 1409--1425.
%
\bibitem{K1}{\sc M.V. Klibanov},  1984, {Inverse problems in the large  
and Carleman bounds}, {\em Differential
Equations}, 20, 755--760.
%
\bibitem{K2}{\sc  M.V. Klibanov}, 1992, Inverse problems and Carleman 
estimates, {\em Inverse Problems}, 8, 575--596.

%
\bibitem{KT}{\sc M.V. Klibanov and A. Timonov}, 2004, {Carleman estimates
for coefficient inverse problems and numerical applications}, 
{\em Inverse and Ill-posed series,VSP, Utrecht}.
%
\bibitem{LTY}{\sc I. Lasiecka, R. Triggiani and P.F. Yao}, 1999,  
{Inverse/observability estimates for second order hyperbolic 
equations with variable coefficients}, {\em J. Math. Anal. Appl.}, 235, 
13--57.
%
\bibitem{MOR}{\sc A. Mercado, A. Osses and L. Rosier}, 2008, {Inverse 
problems for the Schr\"{o}dinger equations via Carleman inequalities 
with degenerate weights}, {\em Inverse Problems}, 24, 015017.
%
\bibitem{PY}{\sc J.P. Puel and M. Yamamoto}, 1997,  {Generic well-posedness 
in a multidimensional hyperbolic inverse problem}, {\em J. Inverse Ill-Posed 
Probl.}, 1, 53--83.

\end{thebibliography}
\end{document}